\documentclass[12pt,leqno]{amsart} 
\newtheorem{thm}{Theorem}[section]
\newtheorem{defi}{Definition}[section]

\newtheorem{cor}{Corollary}[section]
\newtheorem{pr}{Proposition}[section]
\newtheorem{rem}{Remark}[section]
\newcommand{\be}{\begin{equation}}
\newcommand{\ee}{\end{equation}}
\newcommand{\bea}{\begin{eqnarray}}
\newcommand{\eea}{\end{eqnarray}}
\newcommand{\beb}{\begin{eqnarray*}}
\newcommand{\eeb}{\end{eqnarray*}}
\usepackage{amssymb,amsfonts,amsthm,setspace,indentfirst}
\numberwithin{equation}{section}
\begin{document}
%
\title[On generalized Roter type manifolds]{\bf{On generalized Roter type manifolds}}
\author[Absos Ali Shaikh and Haradhan Kundu]{Absos Ali Shaikh and Haradhan Kundu}
\date{}
\address{\noindent\newline Department of Mathematics,\newline University of 
Burdwan, Golapbag,\newline Burdwan-713104,\newline West Bengal, India}
\email{aask2003@yahoo.co.in}
\email{kundu.haradhan@gmail.com}
\dedicatory{Dedicated to Professor Witold Roter on his eighty-second birthday}
\begin{abstract}
The main object of the present paper is to study the geometric properties of a generalized Roter type semi-Riemannian manifold, which arose in the way of generalization to find the form of the Riemann-Christoffel curvature tensor $R$. Again for a particular curvature restriction on $R$ and the Ricci tensor $S$ there arise two structures, e. g., local symmetry ($\nabla R = 0$) and Ricci symmetry ($\nabla S = 0$); semisymmetry($R\cdot R =0$) and Ricci semisymmetry ($R\cdot S =0$) etc. In differential geometry there is a very natural question about the equivalency of these two structures. In this context it is shown that generalized Roter type condition is a sufficient condition for various important second order restrictions. Some generalizations of Einstein manifolds are also presented here. Finally the proper existence of both type of manifolds are ensured by some suitable examples.
\end{abstract}


\maketitle
%
\section{\bf{Introduction}}
In differential geometry, many geometric structures are formed by some curvature restrictions imposed on the Riemann-Christoffel curvature tensor $R$, e.g., locally symmetric manifold, which is formed by $\nabla R = 0$, where $\nabla$ denotes the covariant derivative. Taking a particular restriction and replace the curvature tensor $R$ to a $(0,4)$-tensor, we get different structures. Recently, in \cite{SK14} the present authors studied the equivalency of these structures and also obtained the classification of various curvature restrictions. Now for a certain restriction, we get another geometric structure by imposing the restriction on the Ricci tensor $S$, e.g., $\nabla R = 0$ gives rise to $\nabla S = 0$ (i.e., Ricci symmetric manifold). Consequently a natural question arises about the equivalency of the two geometric structures formed by a particular curvature restriction imposed on $R$ and $S$ respectively. In this context we mention some well known problems:\\
(i) What are the conditions for equivalency of local symmetry and Ricci symmetry?\\
(ii) Are the restrictions of semisymmetry and Ricci semisymmetry equivalent for the hypersurfaces in Euclidean spaces? (P. J. Ryan Problem \cite{Ryan72})\\
(iii) What are the conditions for equivalency of semisymmetry and Ricci semisymmetry? \\
(iv) What are the conditions for equivalency of pseudosymmetry and Ricci pseudosymmetry?\\
It is obvious that if the corresponding restriction operator and contraction commute then the structure due to $R$ implies the structure due to $S$ but not conversely, in general.\\
\indent To answer the above questions, in this paper we study the generalized Roter type manifold (briefly, $GRT_n$) and show that if the corresponding restriction operator is linear over $C^{\infty}(M)$, the ring of all smooth functions on $M$, and commutes with contraction then the two structures formed by this restriction imposing on $R$ and $S$ respectively are equivalent on a $GRT_n$.  We note that in \cite{ACDE98} Arslan et. al. and in \cite{DHS99} Deszcz et. al. presented some sufficient conditions for the above problems (iii) and (iv), which are more complicated than generalized Roter type condition. The papers \cite{ACDE98} and \cite{DHS99} are totally devoted to solve the problems (iii) and (iv) respectively but here we also prove these in other way with the help of our sufficient condition.\\
\indent The notion of $GRT_n$ arises to express $R$ as a simple algebraic combination of some lower order tensors. Consider the following forms of $R$: 
$$(a) R = I g\wedge g,$$
$$(b) R = J_1 g\wedge g + J_2 g\wedge S,$$
$$(c) R =N_1 g\wedge g + N_2 g\wedge S + N_3 S\wedge S,$$
$$(d) R= L_1 g\wedge g + L_2 g\wedge S + L_3 S\wedge S + L_4 g\wedge S^2 + L_5 S\wedge S^2 + L_6 S^2\wedge S^2,$$
where $I,J_i, N_i, L_i$'s are all belong to $C^{\infty}(M)$, and $S^2$ is the Ricci tensor of level 2, and is defined as $S^2(X,Y) = S(\mathcal S X, Y)$, $\mathcal S$ is the Ricci operator and $X,Y\in \chi(M)$, the Lie algebra of all smooth vector fields on $M$.
We note that $I$, $J_1$ and $J_2$ can be determined as $\frac{\kappa}{2 n(n-1)}$, $-\frac{\kappa}{2(n-1)(n-2)}$ and $\frac{1}{n-2}$ respectively, $\kappa$ being the scalar curvature of $M$. A manifold  satisfying the conditions (a), (b), (c) and (d) are respectively called manifold of constant curvature, manifold with vanishing conformal curvature tensor (more precisely conformally flat for $n>3$), Roter type manifold (briefly, $RT_n$) and $GRT_n$. Thus from the definition we have a path of generalization such that corresponding classes are in an inclusion form given as follows:
\begin{center}
\fbox{constant curvature} $\subset$ \fbox{conformally flat} $\subset$ \fbox{Roter type} $\subset$ \fbox{generalized Roter type}.
\end{center}
We mention that a manifold which is not $RT_n$ but satisfying some generalized Roter type condition, was already investigated in \cite{Sawi06} and very recently in \cite{DGJP-TZ13} and \cite{DHJKS13} but the name ``generalized Roter type'' was first used  in \cite{SDHJK15}. Recently in \cite{SKgrtw} the authors studied the characterization of a warped product manifold satisfying some generalized Roter type condition.\\
\indent There are some physical points of view to study this generalized structure. We see that conformally flat manifold has a great importance in general theory of relativity and cosmology. Various spacetimes are conformally flat but some of them are not so, but we need to know their curvature form. Many important spacetimes are neither conformally flat nor $RT_n$ but $GRT_n$. For example, interior black hole spacetime \cite{DHKS} is not conformally flat but a $RT_n$. Again in \cite{DHJKS13} Deszcz et. al. presented some generalization of G\"{o}del metric and showed that this is not $RT_n$ but $GRT_n$ and thus as special cases G\"{o}del metric and Som-Raychaudhuri solution of Einstein's field equation are not $RT_n$ but they are $GRT_n$.\\
\indent The paper is organized in the following way: Section 2 is concerned with preliminaries and section 3 deals with the definitions of various geometric structures and we present some generalizations of Einstein condition. Section 4 is devoted to the main results. Finally, in section 5 we present two examples which ensure the properness of the path of generalization to get the form of curvature tensor $R$ and the path of generalization of Einstein condition.
\section{\bf{Preliminaries}}
In this paper we consider all the manifold to be smooth connected semi-Riemannian and of dimension $n$, $n\ge 3$. Let us consider a manifold $M$ equipped with a semi-Riemannian metric $g$ and $\nabla$, $R$, $S$ and $\kappa$ be the corresponding Levi-Civita connection, Riemann-Christoffel curvature tensor, Ricci tensor and scalar curvature respectively.
Now for two $(0,2)$-tensors $A$ and $E$, their Kulkarni-Nomizu product (\cite{DG02}, \cite{DGHS98}, \cite{DH03}, \cite{Glog02}) $A\wedge E$ is given by
\bea\label{eq2.2}
(A \wedge E)(X_1,X_2,Y_1,Y_2)&=&A(X_1,Y_2)E(X_2,Y_1) + A(X_2,Y_1)E(X_1,Y_2)\\\nonumber
&-&A(X_1,Y_1)E(X_2,Y_2) - A(X_2,Y_2)E(X_1,Y_1),
\eea
where $X_1, X_2, Y_1, Y_2\in \chi(M)$. Throughout the paper we consider $X, Y, X_i, Y_i \in \chi(M)$, $i = 1,2, \cdots $.\\
Again for a $(0,2)$-tensor $A$ and a $(0,k)$-tensor $T$ we define their generalized Kulkarni-Nomizu product (\cite{ADEHM14}, \cite{DG02a}) as a $(0,k+2)$-tensor $A\wedge T$ which is given by
\bea\label{eq2.1}\nonumber
(A \wedge T)(X_1,X_2,Y_1,Y_2,\cdots,Y_k)&=&A(X_1,Y_2)T(X_2,Y_1,\cdots,Y_k) + A(X_2,Y_1)T(X_1,Y_2,\cdots,Y_k)\\
&-&A(X_1,Y_1)T(X_2,Y_2,\cdots,Y_k) - A(X_2,Y_2)T(X_1,Y_1,\cdots,Y_k).
\eea
For a symmetric $(0,2)$-tensor $A$ we get an endomorphism $\mathcal A$ called the corresponding endomorphism operator defined as
$$g(\mathcal AX,Y) = A(X,Y).$$
Then we can define $k$-th level tensor of $A$, say $A^k$ of same order with corresponding endomorphism operator $\mathcal A^k$ given below:
$$A^k(X,Y) = g(\mathcal A^k X,Y) = A(\mathcal A^{k-1}X,Y).$$
Thus, in particular, we get the second, third and fourth level Ricci tensor $S^2$, $S^3$, $S^4$ respectively given by
$$S^2(X,Y) = S(\mathcal SX,Y), \ \ \ S(X,Y) = g(\mathcal SX,Y),$$
$$S^3(X,Y) = S(\mathcal S^{2}X,Y), \ \ \ S^2(X,Y) = g(\mathcal S^2 X,Y),$$
$$S^4(X,Y) = S(\mathcal S^{3}X,Y), \ \ \ S^3(X,Y) = g(\mathcal S^3 X,Y).$$
A tensor $D$ of type $(1,3)$ on $M$ is said to be generalized curvature tensor (\cite{DG02}, \cite{DGHS98}, \cite{DH03}), if
\beb
&(i)&D(X_1,X_2)X_3+D(X_2,X_3)X_1+D(X_3,X_1)X_2=0,\\
&(ii)&D(X_1,X_2)X_3+D(X_2,X_1)X_3=0,\\
&(iii)&D(X_1,X_2,X_3,X_4)=D(X_3,X_4,X_1,X_2),
\eeb
where $D(X_1,X_2,X_3,X_4)=g(D(X_1,X_2)X_3,X_4)$, for all $X_1,X_2,$ $X_3,X_4$. Here we use the same symbol $D$ for the generalized curvature tensor of type $(1,3)$ and $(0,4)$. 
Moreover if $D$ satisfies the second Bianchi like identity i.e.,
$$(\nabla_{X_1}D)(X_2,X_3)X_4+(\nabla_{X_2}D)(X_3,X_1)X_4+(\nabla_{X_3}D)(X_1,X_2)X_4=0,$$
then $D$ is called a proper generalized curvature tensor. 
Some most useful generalized curvature tensors are Gaussian curvature tensor $G$, Weyl conformal curvature tensor $C$, concircular curvature tensor $W$ and conharmonic curvature tensor $K$, which are respectively given by
\begin{eqnarray*}
G &=& \frac{1}{2}(g \wedge g),\\
C &=& R -\frac{1}{(n-2)}(g \wedge S)+\frac{\kappa}{2(n-1)(n-2)}(g \wedge g),\\
W &=& R - \frac{\kappa}{2 n(n-1)}(g \wedge g),\\
K &=& R - \frac{1}{(n-2)}(g \wedge S).
\end{eqnarray*}
We note that the Gaussian curvature tensor is a proper generalized curvature tensor.\\
We can easily operate an endomorphism $\mathcal H$ over $\chi(M)$, on a $(0,k)$-tensor $T$, $k\geq 1$, and get the tensor $\mathcal H T$, given by (\cite{DG02}, \cite{DGHS98}, \cite{DH03})
\beb\label{hdot}
(\mathcal{H} T)(X_1,X_2,\cdots,X_k) = -T(\mathcal{H}X_1,X_2,\cdots,X_k) - \cdots - T(X_1,X_2,\cdots,\mathcal{H}X_k).
\eeb
We note that the operation of $\mathcal H$ on a scalar is zero.\\
Now for a $(0,4)$ tensor $D$ and given two vector fields $X,Y\in\chi(M)$ one can define an endomorphism $\mathcal{D}(X,Y)$ by
$$\mathcal{D}(X,Y)(X_1)=D(X,Y)X_1, \ \mbox{for all $X_1\in\chi(M)$}.$$
Again if $X,Y\in\chi(M)$ then for a $(0,2)$-tensor $A$ one can define an endomorphism $X \wedge_A Y$, by
$$(X \wedge_A Y)X_1 = A(Y,X_1)X - A(X,X_1)Y, \ \mbox{for all $X_1 \in \chi(M) $}.$$
Let $\mathcal T^r_k(M)$ be the space of all smooth tensor fields of type $(r,k)$ on $M$, 
$r,k\in \mathbb N\cup\left\{0\right\}$. Then for $T\in \mathcal T^0_k(M)$, $k\geq 2$, 
and a generalized curvature tensor $D$ one can define a $(0,k+2)$ 
tensor $D\cdot T$ given by (\cite{DDVY94}, \cite{DG02}, \cite{DGHS98}, \cite{DH03}, \cite{SK14})
\beb\label{rdot}
&&D\cdot T(X_1,X_2, \ldots ,X_k;X,Y) = (\mathcal{D}(X,Y)\cdot T)(X_1,X_2, \ldots ,X_k)\\\nonumber
&&= -T(\mathcal{D}(X,Y)X_1,X_2, \ldots ,X_k) - \ldots - T(X_1,X_2, \ldots ,\mathcal{D}(X,Y)X_k),
\eeb
and for a $(0,2)$-tensor $A$ one can define a $(0,k+2)$-tensor $Q(A,T)$ as (\cite{DDVY94}, \cite{DG02}, \cite{DGHS98}, \cite{DH03}, \cite{SK14}, \cite{Tach74})
\beb\label{qgr}
&&Q(A,T)(X_1,X_2, \ldots ,X_k;X,Y) = ((X \wedge_A Y)\cdot T)(X_1,X_2, \ldots ,X_k)\\\nonumber
&&= -T((X \wedge_A Y)X_1,X_2, \ldots ,X_k) - \ldots - T(X_1,X_2, \ldots ,(X \wedge_A Y)X_k).
\eeb
To define the restriction of pseudosymmetric manifold in sense of Chaki \cite{Chak87} as an operator, recently Shaikh and Kundu \cite{SK14} defined an endomorphism $\mu_{_X}$ for an $1$-form $\mu$ and a vector field $X$ on $M$ as
$$\mu_{_X}(X_1) = \mu(X_1)X, \ \mbox{for all $X_1\in \chi(M).$}$$
Then we can operate $\mu_{_X}$ on a $(0,k)$-tensor field $T$ as follows:
\beb\label{pidot}
&&(\mu_{_X} \cdot T)(X_1,X_2, \ldots ,X_k)\\\nonumber
&&= -T(\mu_{_X}(X_1),X_2, \ldots ,X_k) - \ldots - T(X_1,X_2, \ldots ,\mu_{_X}(X_k)),\\\nonumber
&&= -\mu(X_1)T(X,X_2, \ldots ,X_k) -\mu(X_2)T(X_1,X, \ldots ,X_k)- \ldots -\mu(X_k)T(X_1,X_2, \ldots ,X),
\eeb
for all $X, X_i \in \chi(M)$.
\begin{pr}\label{pr2.1}
Let $\mathcal H$ be an endomorphism over $\chi(M)$ and $A$, $E$ be two $(0,2)$-tensors. Then
$$\mathcal H (A\wedge E) = (A\wedge \mathcal H E) + (E\wedge \mathcal H A).$$
\end{pr}
\textbf{Proof:} The proof can be easily follows from the definitions.
\section{\bf{Various geometric structures}}
In this section we define various necessary geometric structures by means of curvature restrictions.
\begin{defi}$($\cite{Cart26}, \cite{SK14}, \cite{Szab82}$)$ 
Let $T$ be a $(0,k)$ tensor and $D$ be a $(0,4)$ tensor on $M$. Then the semi-Riemannian manifold $M$ is said to be $T$-semisymmetric type with respect to $D$ or simply $T$-semisymmetric type if $D\cdot T = 0$.
\end{defi}
In particular, for $D= R$ and if $T= R$ (resp., $S$, $P$, $C$, $W$, $K$) then the manifold is called semisymmetric 
(resp., Ricci semisymmetric, projectively semisymmetric, conformally semisymmetric, concircularly semisymmetric, conharmonically semisymmetric).
\begin{defi}$($\cite{AD83}, \cite{Desz87}, \cite{Desz92}, \cite{SK14}$)$ 
Let $T$ be $(0,k)$ tensor and $D_1, D_2,\cdots D_r$ are some $(0,4)$ tensors on $M$ ($r>1$). Then the semi-Riemannian manifold $M$ is said to be Deszcz $T$-pseudosymmetric type 
if the tensors $D_1 \cdot T$, $D_2 \cdot T$, $\cdots$, $D_r \cdot T$ are linearly dependent.
\end{defi}
In particular, if $r =2$, $D_1 = R$, $D_2 = G$ and $T= R$ (resp., $S$, $P$, $C$, $W$, $K$), then the manifold is called Deszcz pseudosymmetric (resp., Ricci pseudosymmetric, projectively pseudosymmetric, conformally pseudosymmetric, concircularly pseudosymmetric, conharmonically pseudosymmetric). Especially, if $r =2$, $D_1 = C$, $D_2 = G$ and $T =C$, then $M$ is called manifold of pseudosymmetric Weyl tensor. Again if $r =2$, $D_1 = R$, $\mathcal D_2 (X,Y) = X\wedge_S Y$ and $T =R$, then the manifold is called Ricci generalized pseudosymmetric. Thus we can say that a manifold is Deszcz pseudosymmetric, Ricci pseudosymmetric (\cite{ACDE98}, \cite{Desz92}), manifold of pseudosymmetric Weyl tensor (\cite{ACDE98}, \cite{Desz92}) and Ricci generalized pseudosymmetric (\cite{DD91}, \cite{DD91a}, \cite{Desz92}) respectively if and only if
$$R\cdot R = L_R Q(g,R) \ \ \mbox{holds on $U_R = \{x\in M : R-\frac{\kappa}{n(n-1)}G \ne 0$ at $x$\}},$$
$$R\cdot S = L_R Q(g,S) \ \ \mbox{holds on $U_S = \{x\in M : S-\frac{\kappa}{n}g \ne 0$ at $x$\}}$$
$$C\cdot C = L_C Q(g,R) \ \ \mbox{holds on $U_R = \{x\in M : C \ne 0$ at $x$\} and}$$
$$R\cdot R = L Q(S,R) \ \ \mbox{holds on $U_R = \{x\in M : Q(S,R) \ne 0$ at $x$\}}.$$
%
\begin{defi} $($\cite{Desz03}, \cite{Desz03a}$)$ 
Let $T\in \mathcal T^0_4(M)$ and $A, E \in \mathcal T^0_2(M)$. Then the tensor $W(T, A, E)$ is defined by
$$W(T,A,E) = T - N_1 A\wedge A - N_2 A\wedge E - N_3 E\wedge E,$$
and is called a Roter type tensor with $T$, $A$ and $E$, where $N_1, N_2$ and $N_3$ are associated scalars. If on a manifold, $W(T,A,E) = 0$ for some $N_1, N_2, N_3$, then the manifold is said to be a Roter type with $T, A, E$; and $N_1, N_2, N_3$ are called the associated scalars.
\end{defi}
In particular if $W(R,g,S) = 0$ i.e.,
\be\label{rt}
R = N_1 g\wedge g + N_2 g\wedge S + N_3 S\wedge S,
\ee
for some scalars $N_1, N_2, N_3$, then the manifold is simply called Roter type manifold \cite{Desz03} and write briefly here as $RT_n$. In this case we write $W(R,g,S)$ simply as $W(R)$. For more details about $RT_n$ we refer the reader to see \cite{Desz03a}, \cite{DPS13}, \cite{Glog07} and also references therein.
\begin{defi} 
Let $T\in \mathcal T^0_4(M)$ and $A, E, F \in \mathcal T^0_2(M)$. Then the tensor $GW(T, A, E, F)$ is defined by
$$GW(T,A,E,F) = T - L_1 A\wedge A - L_2 A\wedge E - L_3 E\wedge E - L_4 A\wedge F - L_5 E\wedge F - L_6 F\wedge F,$$
and is called a generalized Roter type tensor with $T$, $A$, $E$ and $F$, where $L_i$, $1\le i\le 6$ are associated scalars. If on a manifold, $GW(T,A,E,F) = 0$ for some $L_i$, $1\le i\le 6$, then the manifold is said to be a generalized Roter type with $T, A, E, F$; and $L_i$, $1\le i\le 6$ are called the associated scalars.
\end{defi}
In particular if $GW(R, g, S, S^2) = 0$ i.e.,
\be\label{grt}
R = L_1 g\wedge g + L_2 g\wedge S + L_3 S\wedge S + L_4 g\wedge S^2 + L_5 S\wedge S^2 + L_6 S^2\wedge S^2,
\ee
for some scalars $L_i$, $1\le i\le 6$, then the manifold is called generalized Roter type manifold \cite{SDHJK15} and we briefly write it here as $GRT_n$. In this case we write $GW(R, g, S, S^2)$ simply as $GW(R)$. We note that a $GRT_n$ (resp., $RT_n$) is called proper if $L_6 \ne 0$ (resp., $N_3 \ne 0$) and also called special if one or more than one of their associated scalar(s) are identically zero or some particular value.\\
\indent We have already seen a way of generalization on the basis of curvature form and now we present a parallel way of generalization to get the linear dependency of Ricci tensors of different levels. For this case the first step is Einstein manifold.
\begin{defi} \cite{Bess87} 
A manifold is said to be Einstein if it satisfies the Einstein condition i.e., its Ricci tensor $S$ is some scalar multiple of the metric tensor $g$.
\end{defi}
We note that in this case $S = \frac{r}{n}g$.
\indent Now by contracting the equation (\ref{rt}) and (\ref{grt}), we get respectively 
$a_1 g+a_2 S+a_3 S^2 = 0$ and 
$a_4 g+a_5 S+a_6 S^2+a_7 S^3+a_8 S^4 = 0$,
where $a_1 =2 N_1 (n-1) + N_2 \kappa$, 
$a_2 =N_2 (n-2) + 2 N_3 \kappa-1$, 
$a_3 =-2 N_3$, 
$a_4 = 2 L_1 (n-1) + L_2 \kappa + L_4 \kappa^{(2)}$, 
$a_5 = L_2 (n-2) + 2 L_3 \kappa + L_5 \kappa^{(2)}-1$, 
$a_6 = -2 L_3 + L_4 (n-2) + L_5 \kappa + 2 L_6 \kappa^{(2)}$, 
$a_7 =-2 L_5$ and $a_8 =-2 L_6$, where $\kappa^{(2)}$ is the trace of $S^2$.\\
Thus we get two generalizations of Einstein manifold, namely generalized Einstein condition \cite{Bess87} and are given by
\be\label{gec2}
S^2 + a_1 S + a_2 g = 0 \ \mbox{and}
\ee
\be\label{gec4}
S^4 + a_3 S^3 + a_4 S^2 + a_5 S + a_6 g = 0,
\ee
where $a_1, a_2, a_3, a_4, a_5$ and $a_6$ are some scalars. We can also construct another generalized Einstein condition as
\be\label{gec3}
S^3 + a_7 S^2 + a_8 S + a_9 g = 0,
\ee
where $a_7, a_8$ and $a_9$ are some scalars. According to the maximum level of Ricci tensor we can say the structures given by (\ref{gec2}), (\ref{gec3}) and (\ref{gec4}) as generalized Einstein condition with 2nd, 3rd and 4th Ricci level respectively and are denoted by $Ein(2)$, $Ein(3)$ and $Ein(4)$ respectively. Similarly we can present the Einstein manifold as $Ein(1)$. Finally we get a path of generalization of Einstein manifold by increasing the Ricci level such that their corresponding classes are in an inclusion form given as follows:
\begin{center}
\fbox{Ricci Flat} $\subset$ \fbox{Einstein} $\subset$ \fbox{$Ein(2)$} $\subset$ \fbox{$Ein(3)$} $\subset$ \fbox{$Ein(4)$}.
\end{center}
We note that throughout the paper, by box we denote the corresponding class of manifolds. We can further algebraically extend the way of generalization to get the curvature form and generalized Einstein conditions by using the Ricci tensors of higher level in the definitions. We note that in literature there are many generalization of Einstein manifold sometime also called generalized Einstein metric conditions and for details about these we refer the reader to see \cite{DGHS11} and also references therein.\\
\indent There is an another way of generalization to get the curvature form, which is given by
\begin{center}
\fbox{Flat} $\subset$ \fbox{Constant curvature} $\subset$ \fbox{quasi constant curvature}.
\end{center}
Then the corresponding parallel way to get the form of Ricci tensor is given by
\begin{center}
\fbox{Ricci Flat} $\subset$ \fbox{Einstein} $\subset$ \fbox{quasi Einstein}.
\end{center}
A manifold is said to be of quasi-constant curvature if its Riemann-Christoffel curvature tensor $R$ is given by
$$R = \alpha' G + \beta' g\wedge(\eta'\otimes\eta'),$$
where $\alpha'$, $\beta'$ are some scalar and $\eta'$ is an $1$-form.\\
Again, a manifold is said to be quasi-Einstein if its Ricci tensor is given by
$$S = \alpha g + \beta \eta\otimes\eta,$$
where $\alpha$, $\beta$ are some scalar and $\eta$ is an $1$-form. We note that in this case the Ricci operator can be expressed as
$$\mathcal S = \alpha \mathcal I + \beta \eta_{_V},$$
where $\mathcal I$ is the identity operator and $V$ is the vector field corresponding to $\eta$ and is given by $g(X,V) =\eta(X)$.
\section{\bf{Main Results}}
In this section we discuss various geometric properties of a $GRT_n$ satisfying the generalized Roter type condition (\ref{grt}).
\begin{thm}\label{th4.1}
A $GRT_n$ is a manifold of constant curvature if it is an Einstein manifold.
\end{thm}
\noindent \textbf{Proof:} If the manifold is Einstein, i.e., $S=\frac{\kappa}{n} g$, then
$$S^2 = \left(\frac{\kappa}{n}\right)^2 g.$$
Now putting these values of $S$ and $S^2$ in (\ref{grt}), we get
$$R = \frac{2}{n^4} \left[L_1 n^4 + \kappa \{L_2 n^3 + \kappa ((L_3 + L_4) n^2 + L_5 n \kappa + L_6 \kappa^2)\}\right] G.$$
This completes the proof.
\begin{thm}\label{th4.2}
A $GRT_n$ is a manifold of quasi constant curvature if it is a quasi-Einstein manifold.
\end{thm}
\noindent \textbf{Proof:} Let $M$ be quasi-Einstein with $S = \alpha g + \beta \eta\otimes\eta$. Then
$$S^2 = \bar\alpha g + \bar\beta \eta\otimes\eta,$$
where $\bar\alpha = \alpha^2$ and $\bar\beta = \beta(2\alpha +\beta||\eta||^2)$. Now applying these values in (\ref{grt}), we get
$$R = \alpha' G +\beta' g\wedge(\eta\otimes\eta),$$
where $\alpha' = 2 \left[L_1 + \alpha (L_2 + \alpha (L_3 + L_4 + \alpha (L_5 + L_6 \alpha)))\right]$ and 
$\beta' = \beta[L_2 +\alpha (2(L_3 + L_4) + 3 L_5 \alpha + 4 L_6 \alpha^2) +(L_4 + \alpha (L_5 + 2 L_6 \alpha)) \beta ||\eta||^2]$. 
Thus the manifold becomes of quasi constant curvature.\\
\indent It is clear from definition that $GRT_n$ is a generalization of $RT_n$. We now present a sufficient condition for a $GRT_n$ to be $RT_n$.
\begin{thm}\label{th4.3}
A $GRT_n$ is a $RT_n$ if it satisfies some proper $Ein(2)$ condition.
\end{thm}
\noindent \textbf{Proof:} Let us consider the generalized Einstein condition as
$$S^2 + a_1 S + a_2 g = 0.$$
Then putting the value of $S^2$ in the generalized Roter type condition we get our assertion.\\
\textbf{Note:} Since every proper $RT_n$ is $Ein(2)$, from above we can say that a proper $GRT_n$ is $RT_n$ if and only if it is $Ein(2)$. Similarly from Theorem \ref{th4.1} (resp., \ref{th4.2}) we can say that a proper $GRT_n$ is of constant curvature (resp., quasi-constant curvature) if and only if it is Einstein (resp., quasi-Einstein).
\begin{thm}\label{th4.4}
A proper $RT_n$ satisfying $(\ref{rt})$ is a $GRT_n$ which satisfies $(\ref{grt})$ such that $L_i$'s are related as
$$L_1 = -\frac{b^2 L_6 + 2 b L_4 N_3}{4 N_3^2},$$
$$L_2 = -\frac{a b L_6 + a L_4 N_3 + b L_5 N_3 -2 N_2 N_3^2}{2 N_3^2},$$
$$L_3 = -\frac{a^2 L_6 + 2 a L_5 N_3-4 N_3^3}{4 N_3^2},$$
where $a=(n-2) N_2 + 2 N_3 \kappa$ and $b=2 (n-1) N_1 + N_2 \kappa$.
\end{thm}
\noindent \textbf{Proof:} Since the manifold satisfies (\ref{rt}), contracting  (\ref{rt}) we get
$$a_1 g + a_2 S + a_3 S^2 = 0,$$
where $a_1 =2 N_1 (n-1) + N_2 \kappa$, 
$a_2 =N_2 (n-2) + 2 N_3 \kappa$, 
$a_3 =-2 N_3$.
Now since $N_3\neq 0$, $a_3 \neq 0$, we can evaluate $S^2$ in terms of $S$ and $g$. Now consider the generalized Roter type tensor
$$GW(R)= R - L_1 g\wedge g - L_2 g\wedge S - L_3 S\wedge S - L_4 g\wedge S^2 - L_5 S\wedge S^2 - L_6 S^2\wedge S^2$$
and put the value of $S^2$ in terms of $S$ and $g$. Then equating the reduced form of $GW(R)$ to zero and solving the corresponding system of equations, we get the result.\\
\indent For a proper conformally flat manifold we can get only the form of $R$ but not of $S$. So in generally we can not conclude about the Roter type and generalized Roter type condition satisfying by a conformally flat manifold. We now present a result on the form of Roter type and generalized Roter type conditions for a manifold of constant curvature.
\begin{thm}
A non flat manifold of constant curvature satisfies $(\ref{rt})$ such that $N_i$'s are related as
$$\frac{N_1 n^2 + \kappa (N_2 n + N_3 \kappa)}{n} = \frac{\kappa}{2(n-1)}.$$
It also satisfies $(\ref{grt})$ such that $L_i$'s are related as
$$\frac{L_1 n^4 + \kappa [L_2 n^3 + \kappa \{n^2 (L_3 + L_4)+ L_5 n \kappa + L_6 \kappa^2\}]}{n^3} = \frac{\kappa}{2(n-1)}.$$
\end{thm}
\noindent \textbf{Proof:} The proof of this theorem may be constructed in same style as similar as Theorem \ref{th4.4}.\\
\indent In literature we see that after introduce of a curvature restricted geometric structure due to a restriction on $R$ there parallely arises another geometric structure due to the same restriction on $S$. If the corresponding operator of the restriction and contraction commute then the first structure implies the second. Thus there arises a natural question about the inverse implication or more precisely, the equivalency of these two structures (see the P. J. Ryan problem \cite{Ryan72},  works of Deszcz and his coauthors \cite{ACDE98}, \cite{DHS99} and references therein).  Here we show that the generalized Roter type condition is a sufficient condition for some of them.\\
\indent For this purpose we first consider about the linearity and commutativity with contraction of the defining restriction operator of a geometric structure. If the operator is linear over $R$ only, then it is called 1st type and if it linear over $C^{\infty}(M)$ also, then it is called 2nd type. Again if the operator and contraction are commute then it is called commutative. We note that a restriction is commutative if and only if it gives zero when it applies on $g$ (Lemma 5.1, \cite{SK14}). For details about the classification of such curvature restriction operators we refer the readers to see \cite{SK14}.
\begin{thm}\label{eqv2}
On a $GRT_n$ the geometric structures defined by a commutative 2nd type restriction imposed on $R$ and $S$ are equivalent.
\end{thm}
\noindent \textbf{Proof:} Consider a commutative 2nd type restriction operator $\mathcal L$ and then the corresponding geometric structure due to a tensor $T$ is given by $\mathcal L (T)=0$. Now as $\mathcal L$ is commutative so $\mathcal L (R) = 0$ $\Rightarrow$ $\mathcal L (S) = 0$. So to prove the theorem it is sufficient to show $\mathcal L (S) = 0$ $\Rightarrow$ $\mathcal L (R) = 0$.\\
Now consider $\mathcal L (S) = 0$. Then from Proposition \ref{pr2.1}, $\mathcal L (S\wedge S) = 0$ and again contracting we get $\mathcal L (S^2) = 0$ (since $\mathcal L$ is commutative).\\
Let us now consider the generalized Roter type condition (\ref{grt}). Then as $\mathcal L$ is of 2nd type and commutative, by Proposition \ref{pr2.1}, we have
$$\mathcal L (R) = 0 + L_2 g\wedge \mathcal L(S) + 2 L_3 S\wedge \mathcal L(S) + L_4 g\wedge \mathcal L(S^2) + 
L_5 S^2 \wedge \mathcal L(S) + L_5 S\wedge \mathcal L(S^2) + 2 L_6 S^2\wedge \mathcal L(S^2).$$
Since $\mathcal L (S) = 0$ and $\mathcal L (S^2) = 0$, so from above we get $\mathcal L (R) = 0$. Hence the theorem.\\
\indent Since the curvature restriction of semisymmetric and pseudosymmetric structures are commutative and of 2nd type, from the above theorem we can state the following corollaries:
\begin{cor}
On a $GRT_n$ the semisymmetric and Ricci semisymmetric conditions are equivalent.
\end{cor}
\begin{cor}
On a $GRT_n$ the pseudosymmetric and Ricci pseudosymmetric conditions are equivalent.
\end{cor}
Again since the generalized Roter type condition is generalization of conformally flat and Roter type condition i.e., they are improperly generalized Roter type, from the above theorems we can state the following:
\begin{cor}\label{cor4.3}
On a $RT_n$ or a conformally flat manifold, the semisymmetric (resp., pseudosymmetric) and Ricci semisymmetric (resp., Ricci pseudosymmetric) conditions are equivalent.
\end{cor}
\begin{rem}
We can conclude that in a $GRT_n$ to study any semisymmetric type or pseudosymmetric type condition imposed on $R$, it is sufficient to study the condition on $S$ only, if the corresponding restriction is commutative 2nd type.
\end{rem}
\indent From above we can state the following:
\begin{cor}
On a $GRT_n$, the following implications hold:
$$R\cdot S = 0 \Rightarrow R\cdot R = 0, R\cdot C = 0,$$
$$R\cdot S = L_S Q(g,S) \Rightarrow R\cdot R = L_S Q(g,R), R\cdot C = L_S Q(g,C),$$
$$C\cdot S = 0 \Rightarrow C\cdot R = 0, C\cdot C = 0,$$
$$C\cdot S = L_S Q(g,S) \Rightarrow C\cdot R = L_S Q(g,R), C\cdot C = L_S Q(g,C),$$
$$W\cdot S = 0 \Rightarrow W\cdot R = 0, W\cdot C = 0,$$
$$W\cdot S = L_S Q(g,S) \Rightarrow W\cdot R = L_S Q(g,R), W\cdot C = L_S Q(g,C),$$
$$K\cdot S = 0 \Rightarrow K\cdot R = 0, K\cdot C = 0$$
$$\mbox{and } K\cdot S = L_S Q(g,S) \Rightarrow K\cdot R = L_S Q(g,R), K\cdot C = L_S Q(g,C).$$
\end{cor}
\indent In \cite{ACDE98} Arslan et. al. showed that in a manifold of special generalized Ricci pseudosymmetry with pseudosymmetric Weyl conformal curvature tensor, the notion of Ricci semisymmetry and semisymmetry are equivalent. Again in \cite{DHS99} Deszcz et. al. showed that in a manifold in which $R\cdot R - Q(S,R)$ and $Q(g,C)$ are linearly dependent together with $R\cdot C$ and $Q(S,C)$ are linearly dependent, the notion of pseudosymmetry and Ricci pseudosymmetry are equivalent. They proved these in a broad way. Here we show that these conditions either give some generalized Roter type condition and thus by our Theorem \ref{eqv2} we can easily conclude the results of equivalency otherwise follows from very simple results. We first state some necessary results.
\begin{pr}\label{pr4.1i}\cite{DD91}
Let $(M,g)$, $dim(M)\ge 3$, be a semi-Riemannian manifold. Let $A$ be a non-zero symmetric $(0,2)$-tensor of $Rank\ge 2$ and $B$ be a generalized curvature tensor. Let at $x\in M$, $Q(A, B) = 0$, then\\
(i) $B$ and $A\wedge A$ are linearly dependent if $A(X,Y) \ne \frac{1}{A(V,V)} A(V,X) A(V,Y)$.\\
(ii) ${\displaystyle\sum_{X,Y,Z}} a(X)\mathcal B(X,Y) = 0$ if $A(X,Y) = \frac{1}{A(V,V)} A(V,X) A(V,Y)$,\\
where $V$ is a vector at $x$ such that $A(V,V)\ne 0$.
\end{pr}
\begin{pr}\label{pr4.1ii}\cite{DD91}
Let $(M,g)$, $dim(M)\ge 3$, be a semi-Riemannian manifold. If at a point $x$ in M,
$$S = \mu g+\rho a\otimes a \  \mbox{and } \ \sum_{X,Y,Z} a(X)\mathcal B(X,Y) = 0$$
for some non-zero vector $a$, where $B = R-\gamma G$, $\mu, \rho, \gamma \in \mathbb R$. Then at $x$, we have
$$R\cdot R = \frac{\kappa}{n(n-1)}Q(g,R) \ \ \mbox{and } \ R\cdot R = Q(S,R)- \frac{(n-2)\kappa}{n(n-1)}Q(g,C).$$
\end{pr}
\begin{pr}\label{pr4.3}\cite{DHS99}
Let $(M,g)$, $dim(M)\ge 3$, be a semi-Riemannian Ricci pseudosymmetric manifold ($R\cdot S = L_S Q(g,S)$) in which $R\cdot R - Q(S,R)$ and $Q(g,C)$ are linearly dependent ($R\cdot R - Q(S,R) = L_1 Q(g,C)$) together with $R\cdot C$ and $Q(S,C)$ are linearly dependent ($R\cdot C = L_2 Q(S,C)$). Then if $L_1 \neq 0$ or $L_2 \neq 1$, then $(M,g)$ satisfies the  pseudosymmetry condition $R\cdot R = L_S Q(g,R)$ on $U_S$.
\end{pr}
\noindent \textbf{Proof:} It is given that
$$R\cdot S = L_S Q(g,S),$$
$$R\cdot R - Q(S,R) = L_1 Q(g,C) \ \  \mbox{and}$$
$$R\cdot C = L_2 Q(S,C).$$
Now comparing last two results we get
$$Q(S,R) + L_1 Q(g,C) - \frac{1}{n-1} R\cdot g\wedge S = L_2 Q(S,C)$$
\beb
\Rightarrow \ Q(S,R) &+& L_1 \left[Q(g,R) - \frac{1}{n-1}Q(g,g\wedge S)\right] - \frac{1}{n-1} g\wedge (R\cdot S)\\
&=& L_2 \left[Q(S,R)- \frac{1}{n-1}Q(S,g\wedge S)+\frac{\kappa}{2(n-1)(n-2)} Q(S,g\wedge g)\right]
\eeb
\beb
\Rightarrow \ Q(S,R) &+& L_1 \left[Q(g,R) + \frac{1}{n-1}Q(S,g\wedge g)\right] - \frac{L_S}{n-1} Q(g, g\wedge S)\\
&=& L_2 \left[Q(S,R) + \frac{1}{n-1}Q(g,S\wedge S)+\frac{\kappa}{2(n-1)(n-2)} Q(S,g\wedge g)\right]
\eeb
$$\Rightarrow (1-L_2)Q(S,R) + L_1 Q(g,R) + \left[\frac{L_1 + L_S}{n-1}-\frac{L_2\kappa}{2(n-1)(n-2)}\right]Q(S,g\wedge g) - \frac{L_2}{n-1} Q(g, S\wedge S) = 0$$
Then we can write
$$\alpha_1 Q(S,R) + \alpha_2 Q(g,R) + \alpha_3 Q(S,g\wedge g) +\alpha_4 Q(g, S\wedge S) = 0$$
\beb
\Rightarrow\left\{\begin{array}{l}
Q\left(\alpha_1 S + \alpha_2 g, R + \frac{\alpha_3}{\alpha_1} g\wedge g + \frac{\alpha_4}{\alpha_2} S\wedge S\right) = 0 \ \ \mbox{if $\alpha_1 \ne 0$ and $\alpha_2\ne 0$}\\
Q\left(\alpha_2 g, R - \frac{\alpha_3}{\alpha_2} g\wedge S + \frac{\alpha_4}{\alpha_2} S\wedge S\right) = 0 \ \ \mbox{if $\alpha_1 = 0$ and $\alpha_2\ne 0$}\\
Q\left(\alpha_1 S, R + \frac{\alpha_3}{\alpha_1} g\wedge g - \frac{\alpha_4}{\alpha_1} g\wedge S\right) = 0 \ \ \mbox{if $\alpha_1 \ne 0$ and $\alpha_2= 0$,}
\end{array}\right.
\eeb
where $\alpha_1 = 1-L_2$, $\alpha_1 = L_1$, $\alpha_3 = \frac{L_1 + L_S}{n-1}-\frac{L_2\kappa}{2(n-1)(n-2)}$ and $\alpha_4 = - \frac{L_2}{n-1}$.
Let us now consider the above results as $Q(A,B) = 0$.
Then two cases arise:\\
\textbf{Case 1:} Let at some $x\in M$, $Rank(A)>1$. Then from Proposition \ref{pr4.1i} we have 
$$B = \lambda_1 A\wedge A$$
i.e., $M$ satisfies some Roter type condition at $x$. Then from Theorem \ref{eqv2}, we can say that the manifold satisfies the pseudosymmetry condition at $x$, since $M$ satisfies the Ricci pseudosymmetry condition.\\
\textbf{Case 2:} Let at some $x\in M$, $Rank(A)=1$ and $A\ne 0$. Then there exists some $V\in T_x(M)$ such that $A(V,V) \ne 0$ and we have
$$\alpha_1 S + \alpha_2 g = \frac{1}{\rho} a\otimes a \ \ \mbox{and}$$
$$a(X)\mathcal B(Y,Z)+a(Y)\mathcal B(Z,X)+a(Z)\mathcal B(X,Y) =0,$$
where $a(X)=A(X,V)$ and $\rho = A(V,V)$. We note that since $g$ can not be of rank 1, the fact $\alpha_1 = 0$ does not arise in this case. Thus $B$ is some linear combination of $R$, $G$ and $g\wedge (a\otimes a)$, say
$$B=\beta_1 R + \beta_2 G+\beta_3 g\wedge (a\otimes a).$$
Since $\tilde B = g\wedge (a\otimes a)$ satisfies the condition
$$a(X)\tilde{\mathcal B}(Y,Z)+a(Y)\tilde{\mathcal B}(Z,X)+a(Z)\tilde{\mathcal B}(X,Y) =0,$$
$\bar B = \beta_1 R + \beta_2 G$ satisfies
$$a(X)\bar{\mathcal B}(Y,Z)+a(Y)\bar{\mathcal B}(Z,X)+a(Z)\bar{\mathcal B}(X,Y) =0.$$
Thus by Proposition \ref{pr4.1ii}, we have
$$R\cdot R = \frac{\kappa}{n(n-1)}Q(g,R).$$
Hence at each $x\in M$, the manifold is pseudosymmetric except $\alpha_1 S + \alpha_2 g = 0$ i.e., on the set $\left\{x\in M: S=\frac{\kappa}{n} \ \mbox{at $x$}\right\}$. This completes the proof.\\
\begin{rem}
We note that in \cite{DHS99} Deszcz et. al. proved the above result by considering $L_S = \frac{\kappa}{n}L_2 \neq 0$. They also showed that under this given condition $L_S = \frac{\kappa}{n(n-1)}$ and thus $L_2 = \frac{1}{n-1}$ (see Remark 3.1, \cite{DHS99}). Since here we consider $n\ge 3$, so $L_2 \neq 1$, thus the result presented by Deszcz et. al. in \cite{DHS99} is generalized by the above proposition.
\end{rem}
\begin{pr}\label{pr4.4}\cite{ACDE98}
Let $(M,g)$, $dim(M)\ge 3$, be a semi-Riemannian Ricci semisymmetric manifold satisfies generalized Ricci pseudosymmetry condition $R\cdot R = Q(S,R)$ with pseudosymmetric Weyl conformal curvature tensor. Then $(M,g)$ satisfies semisymmetry condition on $U_S$.
\end{pr}
\noindent \textbf{Proof:} First consider $\kappa = 0$, the result is easily follows from Theorem 4.2 of \cite{ACDE98}. Next consider $\kappa \neq 0$. Now since $M$ is Ricci semisymmetric and $R\cdot R = Q(S,R)$, from Proposition 5.1 of \cite{ACDE98} we have
$$C\cdot C = \frac{n-3}{n-2}R\cdot R + \frac{1}{n-2}\left(\frac{\kappa}{n-1} - \tau\right) Q(g,C),$$
where $\tau = \frac{tr(S^2)}{\kappa}$. Again the manifold is of pseudosymmetric Weyl conformal curvature tensor i.e.,
$$C\cdot C = L Q(g,C)$$
for some scalar L (say). Now comparing these two results together with generalized Ricci pseudosymmetry condition we get
$$\frac{n-3}{n-2}Q(S,R) + \frac{1}{n-2}\left(\frac{\kappa}{n-1} - \tau\right) Q(g,C) = L Q(g,C)$$
$$\Rightarrow \ \frac{n-3}{n-2}Q(S,R) + \frac{1}{n-2}\left(\frac{\kappa}{n-1} - \tau -L\right) \left[Q(g,R) - \frac{1}{n-1}Q(g,g\wedge S)\right] = 0$$
$$\Rightarrow \ \frac{n-3}{n-2}Q(S,R) + \frac{1}{n-2}\left(\frac{\kappa}{n-1} - \tau -L\right) \left[Q(g,R) - \frac{1}{n-1}Q(S,G)\right] = 0$$
$$\Rightarrow \ Q(\alpha_1 S + \alpha_2 g,R+\alpha_3 G)=0,$$
where $\alpha_1 = \frac{n-3}{n-2}$, $\alpha_2 = \frac{1}{n-2}\left(\frac{\kappa}{n-1} - \tau -L\right)$ and $\alpha_3 = \frac{1}{(n-1)(n-3)}\left(\frac{\kappa}{n-1} - \tau -L\right)$. If $n=3$, $C=0$ and we get our assertion from Corollary \ref{cor4.3} and if $n>3$, $\alpha_1\neq 0$, then as similar to the proof of Proposition \ref{pr4.3}, we get our assertion.\\
\indent We now show some special generalized Roter type conditions as a sufficient condition for equivalency of another two structures, namely, locally symmetric ($\nabla R = 0$) and Ricci symmetric ($\nabla S = 0$).
\begin{thm}\label{eqv1}
On a $GRT_n$ locally symmetry and Ricci symmetry are equivalent if the associated scalars of $GRT_n$ are constant or constant multiple of $\kappa$.
\end{thm}
\noindent \textbf{Proof:} The proof is similar to the proof of the Theorem \ref{eqv2}.
\begin{cor}\label{cor4.4}
On a conformally flat manifold locally symmetry and Ricci symmetry are equivalent.
\end{cor}
\section{\bf{Examples}}
In this section we present a metric to ensure the existence of $GRT_n$ and the manifold with various generalized Einstein conditions.\\
\textbf{Example 1:} Let $M_1$ be a $5$-dimensional connected semi-Riemannian manifold endowed with the semi-Riemannian metric
\be\label{met1}
ds^2 = f\left[(dx^1)^2+(dx^2)^2+(dx^3)^2+(dx^4)^2+ h (dx^5)^2\right],
\ee
where $f$ is a smooth function of $x^1$ and $h$ is a smooth function of $x^1$ and $x^2$. Then its non-zero Riemann-Christoffel curvature tensor and Ricci tensor components (upto symmetry) are the following:

$$R_{1212}= R_{1313}= R_{1414}=\frac{\left(f'\right)^2-f f''}{2 f}, \ \ R_{2323}= R_{2424}= R_{3434}=-\frac{\left(f'\right)^2}{4 f},$$

$$R_{1515}=\frac{1}{4} \left(-2 h f''-h_1 f'+\frac{2 h \left(f'\right)^2}{f}+\frac{f h_1^2}{h}-2 f h_{11}\right), \ \ R_{1525}=\frac{1}{4} f \left(\frac{h_1 h_2}{h}-2 h_{12}\right),$$

$$R_{2525}=\frac{1}{4} f \left(-\frac{f' \left(h f'+f h_1\right)}{f^2}+\frac{h_2^2}{h}-2 h_{22}\right), \ \ R_{3535}= R_{4545}=-\frac{f' \left(h f'+f h_1\right)}{4 f},$$

$$S_{11}=\frac{-f^2 h_1^2+f h h_1 f'+2 h \left[4 f h f''+f^2 h_{11}-4 h \left(f'\right)^2\right]}{4 f^2 h^2}, \ \ 
S_{12}=-\frac{h_1 h_2-2 h h_{12}}{4 h^2},$$

$$S_{22}=\frac{h \left[2 f h f''+2 f^2 h_{22}+f h_1 f'+h \left(f'\right)^2\right]-f^2 h_2^2}{4 f^2 h^2}, \ \ 
S_{33}= S_{44}=\frac{2 f h f''+h \left(f'\right)^2+f h_1 f'}{4 f^2 h},$$

$$S_{55}=\frac{-f^2 \left(h_1^2+h_2^2\right)+2 f h \left[2 h_1 f'+f \left(h_{11}+h_{22}\right)\right]+h^2 \left[2 f f''+\left(f'\right)^2\right]}{4 f^2 h},$$
where $f' = \frac{df}{dx^1}$, $f'' = \frac{d}{dx^1}(\frac{df}{dx^1})$, $h_1 = \frac{\partial h}{\partial x^1}$, $h_2 = \frac{\partial h}{\partial x^2}$, $h_{11} = \frac{\partial }{\partial x^1}(\frac{\partial h}{\partial x^1})$, $h_{12} = \frac{\partial }{\partial x^1}(\frac{\partial h}{\partial x^2})$ and $h_{22} = \frac{\partial }{\partial x^2}(\frac{\partial h}{\partial x^2})$.

Then by stateforward calculations we can calculate the components of $S^2$, $S^3$, $S^4$, $g\wedge g$, $g\wedge S$, $S\wedge S$, $g\wedge S^2$, $S\wedge S^2$ and $S^2\wedge S^2$. Now we discuss the results by taking restrictions on the functions $f$ and $h$.\\
(i) Let $h_1\ne 0$, $h_2\ne 0$ and $f$ is non-constant. Then the manifold is $Ein(4)$ but not $GRT_5$.\\
(ii) Now let $f$ and $h$ both are non-constant function of $x^1$ only. Then $M_1$ is $Ein(3)$ and $GRT_5$.\\
(iii) If $f$ is non-constant such that $15 (f')^3-18 f f' f''+4f^2 f'' \ne 0$ and $h =\frac{c_1 (f')^2}{f^3}$ , $c_1$ is any arbitrary constant, then the manifold is proper $RT_5$ and thus also $Ein(2)$.\\
(iv) If $f$ is non-constant such that $15 (f')^3-18 f f' f''+4f^2 f'' = 0$ and $h =\frac{c_1 (f')^2}{f^3}$, then the manifold is proper conformally flat, where one of $21 (f')^3-22 f f' f''+4f^2 f''$, $33 (f')^3-46 f f' f''+12f^2 f''$ and $-9 (f')^3-2 f f' f''+4f^2 f''$ are not zero all together.\\
(v) If $f$ is non-constant and $h =\frac{c_1 (f')^2}{f^3}$ and also $21 (f')^3-22 f f' f''+4f^2 f'' = 33 (f')^3-46 f f' f''+12f^2 f'' = -9 (f')^3-2 f f' f''+4f^2 f''=0$. Then $M_1$ become of non-flat constant curvature.\\
(vi) Finally if both $f$ and $h$ are constant then the manifold becomes flat.\\
\textbf{Note 1:} From the above example we see that the steps of paths of generalization
\begin{center}
\fbox{Flat} $\subset$ \fbox{constant curvature} $\subset$ \fbox{conformally flat} $\subset$ \fbox{proper $RT_n$} $\subset$ \fbox{proper $GRT_n$}
\end{center}
\begin{center}
and \fbox{Ricci Flat} $\subset$ \fbox{Einstein} $\subset$ \fbox{$Ein(2)$} $\subset$ \fbox{$Ein(3)$} $\subset$ \fbox{$Ein(4)$}.
\end{center}
are all proper.\\
\indent We now present a metric to ensure that the generalized Einstein conditions $Ein(2)$ and $Ein(4)$ are proper generalization of $RT_n$ and $GRT_n$ respectively.\\
\textbf{Example 2:} Let $M_2$ be an open connected subset of $\mathbb R^6$ where $x^5 > 0$, endowed with the semi-Riemannian metric
\be\label{met2}
ds^2 = (dx^1)^2 + e^{x^1}(dx^2)^2 + e^{x^1}(dx^3)^2 + (dx^4)^2 + e^{x^4}(dx^5)^2 + e^{x^4}(x^5+1)^2(dx^6)^2.
\ee
Then by easy calculation we get the non-zero components of Riemann-Christoffel curvature tensor and Ricci tensor (upto symmetry) as:
$$R_{1212}= R_{1313}=-\frac{e^{x^1}}{4}, \ \ R_{2323}=-\frac{1}{4} e^{2 x^1}, \ \ R_{4545}=-\frac{e^{x^4}}{4}, \ \ 
e^{x^4} R_{4646}= R_{5656}=-\frac{1}{4} e^{2 x^4} \left(x^5+1\right)^2$$

$$S_{11}= S_{44}=\frac{1}{2}, \ \ S_{22}= S_{33}=\frac{e^{x^1}}{2}, \ \ S_{55}=\frac{e^{x^4}}{2}, \ \ S_{66}=\frac{1}{2} e^{x^4} \left(x^5+1\right)^2$$

Then we can easily shown that the manifold is Einstein and also satisfies the other generalized Einstein conditions as
$$a_0 g_{_2} + a_1 S_{_2} + a_2 S^2_{_2} = 0,$$
$$a_3 g_{_2} + a_4 S_{_2} + a_5 S^2_{_2} + a_6 S^3_{_2} = 0,$$
$$a_7 g_{_2} + a_8 S_{_2} + a_9 S^2_{_2} + a_{10} S^3_{_2} + a_{11} S^4_{_2} = 0,$$
where 
$a_2 = -2 (2 a_0 + a_1)$, 
$a_6 = -2 (4 a_3 + 2 a_4 + a_5)$, 
$a_{11} = -16 a_7 - 8 a_8 - 4 a_9 - 2 a_{10}$
and $a_0$, $a_1$, $a_3$, $a_4$, $a_5$, $a_7$, $a_8$, $a_9$, $a_{10}$ are arbitrary scalars.\\
We can easily check that the manifold is not a $GRT_6$. Although we get some dependency of $g\wedge g$, $g\wedge S$, $S\wedge S$, $g\wedge S^2$, $S\wedge S^2$ and $S^2\wedge S^2$, given below:
$$L_0 g\wedge g + L_1 g\wedge S + L_2 S\wedge S = 0,$$
$$L_3 g\wedge g + L_4 g\wedge S + L_5 S\wedge S + L_6 g\wedge S^2 + L_7 S\wedge S^2 +L_8 S^2\wedge S^2 =0,$$
where
$L_2 = -4 (L_0 + L_1),$ 
$L_8 = -4 (4 L_3 + 4 L_4 + L_5 + 2 L_6 + L_7)$ 
and $L_i$, $i=0,1,3,4,5,6,7$ are arbitrary scalars.\\
\textbf{Note 2:} From the above example we see that the steps of the following paths of generalization
\begin{center}
\fbox{constant curvature} $\subset$ \fbox{Einstein},
\end{center}
\begin{center}
\fbox{proper $RT_n$} $\subset$ \fbox{$Ein(2)$} and
\end{center}
\begin{center}
\fbox{proper $GRT_n$} $\subset$ \fbox{$Ein(4)$}
\end{center}
are all proper.\\
\textbf{Note 3:} From the Note 1 and 2, we can conclude about the properness at each step of the following net of generalization
\begin{center}
\begin{tabular}{ccccccccc}

\fbox{constant curvature} & $\subset$ & \fbox{$C=0$} & $\subset$ & \fbox{proper $RT_n$} & & $\subset$ & & \fbox{proper $GRT_n$}\\
$\cap$ & & & & $\cap$ & & & & $\cap$\\
\fbox{Einstein or $Ein(1)$} & & $\subset$ & & \fbox{$Ein(2)$} & $\subset$ & \fbox{$Ein(3)$} & $\subset$ & \fbox{$Ein(4)$}

\end{tabular}
\end{center}
\section{\bf{Conclusion}}
In this paper we obtain the forms of Roter type and generalized Roter type conditions of a manifold of constant curvature but not for a conformally flat manifold. For conformally flat case the coefficients of $S\wedge S$, $g\wedge S^2$, $S\wedge S^2$, $S^2\wedge S^2$ may or may not be zero. We present a way of generalization to get the curvature form and a way of generalization of Einstein condition.\\
\indent We also showed (Theorem \ref{eqv2}) that on a $GRT_n$ the two geometric structures formed by a commutative 2nd type restriction (i.e., the defining restriction operator commute with contraction and linear over $C^{\infty}(M)$) imposed on $R$ and $S$ respectively are equivalent. Thus to check a manifold to occupy the structure formed by a commutative 2nd type restriction imposed on $R$ we have just to check the corresponding condition for $S$ only if the manifold is a $GRT_n$, and hence the work becomes more easy. Again, if we extend the generalized Roter type condition to higher Ricci level, then we can conclude that the above equivalency also remain.\\
\noindent
\textbf{Acknowledgment:} 
The second named author gratefully acknowledges to CSIR, 
New Delhi (File No. 09/025 (0194)/2010-EMR-I) for the financial assistance. All the algebraic computations of Section 5 are performed by a program in Wolfram Mathematica.


\end{document}